\documentclass[12pt,reqno]{amsart}
\usepackage{amsmath,amssymb,amsthm,graphicx,epstopdf, xcolor}
\usepackage{mathpazo,verbatim, fullpage} 
\usepackage[all]{xy} 
\newcommand{\transp}[2]{#1 \, #2}
\newcommand{\cyclett}[6]{#1 \, #2 \, .\, #3 \, #4 \, .\, #5 \, #6}
\newcommand{\HM}{\mathbb H} 
\newcommand{\FF}{\mathbb F} 
\newcommand{\MM}{\mathbb M} 
\newcommand{\conic}{\mathcal K} 
\newcommand{\ra}{\rightarrow} 
\newcommand{\lra}{\longrightarrow} 
\newcommand{\SG}{S} 
\newcommand{\num}{\textsc{num}}
\newcommand{\ltr}{\textsc{ltr}}
\newcommand{\PL}{{\mathbb P}^1}
\newcommand{\oo}[1]{{{{\rm a}{#1}}}}
\newcommand{\pP}[4]{{\rm #1}#2_{#3#4}}
\newcommand{\kirk}[3]{N{#1}_{#2#3}} 
\newcommand{\pasc}[3]{L{#1}_{#2#3}} 
\newcommand{\hkirk}[4]{N{#1}_{#2#3}^{(#4)}} 
\newcommand{\hpasc}[4]{L{#1}_{#2#3}^{(#4)}} 
\newcommand{\PLN}[4]{N{#1}{#2}.{#3}'{#4}'} 
\newcommand{\SVL}[4]{L{#1}{#2}.{#3}'{#4}'} 
\newcommand{\KKK}[3]{{\mathbb K}[#1,#2#3]}
\newcommand{\PPP}[3]{{\mathbb P}[#1,#2#3]}
\newcommand{\LLL}[4]{\Lambda[#1#2.#3'#4']}
\newcommand{\MMM}[4]{{\mathbb M}[#1#2.#3'#4']}
\newcommand{\Proof}{\emph{Proof.} \;}
\newcommand{\ask}[1]{\textcolor{red}{(ASK) #1}}
\newcommand{\response}[1]{\textcolor{blue}{Response by JC: \; #1}}
\newcommand{\respond}[2]{\textcolor{blue}{Response by #1: \; #2}}
\newtheorem{Theorem}{Theorem}[section]
\newtheorem{Lemma}[Theorem]{Lemma}
\newtheorem{Proposition}[Theorem]{Proposition}

\begin{document} 
\title{Absolute Projectivities \\ in \\ Pascal's Multimysticum} 
\author{Jaydeep Chipalkatti and Alex Ryba} 
\maketitle


\medskip 

\parbox{16.5cm}{ \small
{\sc Abstract:} Let $\conic$ denote a nonsingular conic in the complex
projective plane. Given six distinct points
$a,b,c,d,e,f$ on $\conic$, Pascal's theorem says that the intersection
points $ab \cap de, bc \cap ef, cd \cap af$ are collinear. The line
containing them is called the Pascal line of the hexagon $abcdef$.
We get sixty such lines by permuting the vertices. These
lines satisfy some incidence theorems, which lead to an 
extremely rich combinatorial structure called the \emph{hexagrammum
  mysticum}. In 1877, Veronese discovered a procedure to create an infinite
sequence of such structures via successive mutations; this sequence is
called the \emph{multimysticum}. In this paper we prove a strong
rigidity theorem which shows that the multimysticum contains $300$ projective ranges that are
absolutely invariant. The two
interspersed sequences of cross-ratios which encode this
invariance turn out to be the alternate convergents in the continued fraction
expansions of $1 \pm 1/\sqrt{3}$.} 

\medskip 



\medskip 

\tableofcontents 

\section{Introduction} 

\subsection{} Pascal's theorem is one of the oldest and loveliest theorems in
classical projective geometry. Given a hexagon inscribed in a 
nonsingular plane conic, it defines the {\em Pascal
line} of the hexagon. We get a collection of sixty such lines by 
permuting the vertices of the hexagon in all possible ways. These lines satisfy
concurrency theorems that define a collection of new points. The new points
in turn satisfy collinearity theorems, eventually leading to an immensely intricate
geometric structure called the \emph{Hexagrammum Mysticum}. (All of this will be  
explained in much greater detail later -- see Section~\ref{section.hm}
below.) The structure, henceforth denoted by $\HM$, consists of $95$ points and
95 lines in the plane. 

For our purposes, it will be convenient to see
$\HM$ as having a {\em fixed} substructure $\FF$ consisting of $35$ points
and $35$ lines, and a {\em mutable} substructure $\MM$
consisting of $60$ points and $60$ lines. The original Pascal lines
are the lines of $\MM$. Thus, at a broad level of conceptualisation,  we have a decomposition 
\[ \HM = \FF + \MM. \] 
There are highly regular incidence relations between the geometric
elements which lie entirely within either 
$\FF$ or $\MM$, and also those which go `across' $\FF$ and $\MM$. 
These incidences were discovered in the nineteenth century
through the combined efforts of several
mathematicians; specifically by Steiner, Pl\"ucker, Kirkman, Cayley and Salmon
(see Sections~\ref{section.kirkman.nodes}--\ref{section.PluckerSalmon} below). 

\subsection{} 
For reasons which will be evident shortly, let us write $\MM^{(0)} =
\MM$ and $\HM^{(0)} = \HM$. In 1877, Veronese discovered a
conceptually elegant inductive procedure to define an infinite sequence of structures 
$\MM^{(i)}$ for $i \geqslant 1$. Starting from the
collection $\MM^{(i)}$ of $60$ points and $60$ lines, this procedure
gives a new collection $\MM^{(i+1)}$ of the same number of points and lines. 

Thus, the sequence of structures 
\[ \HM^{(i)} = \FF + \MM^{(i)}, \] 
can be seen as a \emph{fixed} part $\FF$ together with an infinite sequence of
\emph{mutations}  
\begin{equation} 
\MM^{(0)} \lra \MM^{(1)} \lra \MM^{(2)} \lra \dots \lra 
\MM^{(i)} \lra \MM^{(i+1)} \lra \dots \end{equation} 
for $i \geqslant 0$. The geometric elements  
(i.e., points and lines) in $\MM^{(i)}$ are said to be at
\emph{height} $i$.  
Given two distinct indices $i$ and $j$, the lines and
points in $\MM^{(i)}$ are generally different from those in
$\MM^{(j)}$. However, the surprising and beautiful property of
Veronese's construction is that, 

\begin{itemize} 
\item 
the formal incidence relations within the elements of each $\MM^{(i)}$, and 
\item 
the formal incidence relations going across $\FF$ and $\MM^{(i)}$, 
\end{itemize} 
are the same for all $i$. This invariance allows us to construct
four distinct ranges\footnote{
By way of terminology, if $X$ is a variety abstractly isomorphic to $\PL$,
then by a range on $X$ we will mean a sequence of points 
$(x_1, x_2,x_3,\dots)$ on $X$. For us, $X$ will either be the set of
points on a line $\Lambda$ in the plane, or the pencil of lines
passing through a fixed point $Q$ in the plane. We will say that
the range is based upon $\Lambda$ or $Q$, as the case may be. 
Two ranges $(x_1, x_2, \dots)$ and $(y_1, y_2, \dots)$ on $X$ and $Y$ respectively are said to be isomorphic, if there is an isomorphism $X \ra Y$ carrying each $x_i$ to $y_i$.} associated to the infinite sequence 
\begin{equation} 
\HM^{(0)}, \quad \HM^{(1)}, \quad \HM^{(2)}, \quad \HM^{(3)}, \quad \dots 
\label{multimysticum} \end{equation}
They will be defined below in Sections~\ref{section.hm} and~\ref{section.multimysticum}. 
\subsection{} 
Our main theorem says that all four of these ranges are
isomorphic to each other, and in fact each is isomorphic to the range 
\begin{equation} 
\infty, \quad 0, \quad 1, \quad \frac{1}{2}, \quad \frac{3}{2},
\quad \frac{3}{7}, \quad \frac{11}{7}, \quad \frac{11}{26}, \quad
\frac{41}{26}, \quad \frac{41}{97}, \quad \frac{153}{97}, \quad
\frac{153}{362} \dots 
\label{proj.range1} \end{equation} 
on the projective line $\PL$. This sequence is obtained by interlacing the two sequences 
\[  1, \quad  \frac{3}{2}, \quad \frac{11}{7}, \quad 
\frac{41}{26}, \quad \frac{153}{97}, \quad \dots \qquad \text{and} \qquad 
\frac{1}{2}, \quad \frac{3}{7},\quad \frac{11}{26}, \quad
\frac{41}{97}, \quad \frac{153}{362}, \quad \dots \] 
which are the alternate convergents in the continued fraction
expansions of the numbers $1+1/\sqrt{3}$ and $1-1/\sqrt{3}$ respectively. 

Henceforth, (\ref{proj.range1}) will be referred to as the \emph{Veronese
  sequence} 
\begin{equation} 
\infty, \quad 0, \quad \alpha_0, \quad \alpha_1, \quad 
\alpha_2, \quad \alpha_3, \quad \dots 
\label{Veronese.sequence} \end{equation} 
where $\alpha_0 =1$, and 
$\alpha_i + \alpha_{i+1} = 2$ if $i$ is odd, and 
$1/\alpha_i + 1/\alpha_{i+1} = 3$ if $i$ is even. The recursive rule 
depends on the parity of $i$, essentially because the rule which
defines the mutation $\MM^{(i)} \lra \MM^{(i+1)}$ does so (see
Section~\ref{section.multimysticum} below). 

The sequence in (\ref{multimysticum}) is usually called the
\emph{multimysticum}. Our result can be interpreted as saying that the
intrinsic projective structure of the multimysticum is `rigid'; that
is to say, it remains 
unaffected by how we choose the initial six points on the conic. Moreover, the
even and odd subsequences 
\[ \HM^{(0)}, \; \HM^{(2)}, \; \HM^{(4)}, \dots \quad \text{and} \quad 
\HM^{(1)}, \; \HM^{(3)}, \; \HM^{(5)}, \dots \] 
in the multimysticum have well-defined `limiting' positions. 

\subsection{} 
The \emph{hexagrammum mysticum} is one of the gems of
combinatorial geometry. Even so, its study can be rather forbidding
because it involves a large
number of points and lines with intricate incidences between
them. Hence, we have chosen to begin with an elementary overview of 
the whole subject. The main theorem is in several parts; we will be
able to formulate it precisely from Section~\ref{section.kirkman.range}
onwards after the requisite notation is available. 

A comprehensive modern introduction
to the \emph{hexagrammum mysticum} as well as the \emph{multimysticum} may be found in the two
papers by Conway and Ryba~\cite{ConwayRyba1, ConwayRyba2}. More
classical references will be given in Section~\ref{section.classical.references}
below. We refer the reader to~\cite{Coxeter} for the necessary basic notions in
projective geometry. 

\section{The Hexagrammum Mysticum} 
\label{section.hm} 

\subsection{Pascal Lines} 
\begin{figure}
\includegraphics[width=12cm]{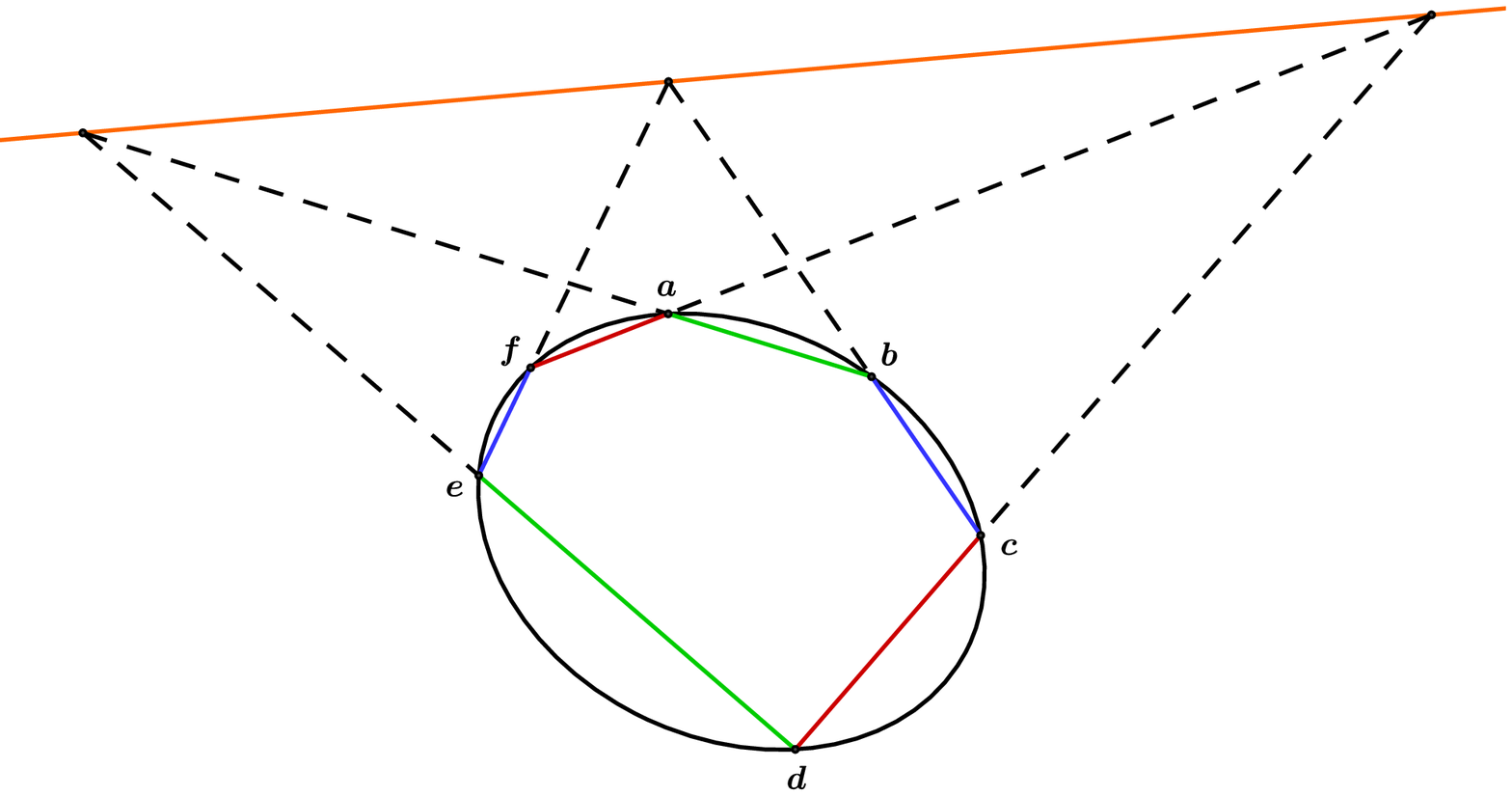}
\caption{\small The Pascal~$L(abcdef)$} 
\label{diagram.pascal.theorem} 
\end{figure}

Let us fix a nonsingular conic $\conic$ in the complex projective
plane, and consider the hexagon formed by six distinct points
$a,b,c,d,e,f$ on $\conic$. Then Pascal's theorem says that the three
intersection points 
\[ ab \cap de, \quad bc \cap ef, \quad cd \cap af, \] 
formed by the three pairs of opposite sides of the hexagon, are
collinear (see Diagram~\ref{diagram.pascal.theorem}). 
The line containing them is called
the Pascal line (or just the Pascal) of the hexagon $abcdef$; we will denote it by 
$L(abcdef)$. If we shift the vertices cyclically or reverse their
order, then the hexagon and its Pascal evidently do not change; hence there are $12$ distinct ways 
\[ L(abcdef) = L(bcdefa) = L(fedcba) \dots \text{etc} \] 
of denoting the same Pascal. However, any essentially different
arrangement of the vertices such as $L(acedbf)$ will \emph{a priori}
correspond to a different Pascal. A moment's reflection will show that
there are $6!/12 = 60$ such arrangements. Henceforth 
we will assume\footnote{We are about to
  define a large number of points and lines with a host of incidence
  relations between them. For some special
  choices of the initial sextuple, some of these
  geometric elements may become undefined or there may be `unwanted'
  incidences between them. However, there is a dense open
  subset of choices $(a,b,c,d,e,f) \in \conic^6$ for which none of
  these pathologies will occur. Henceforth we will always assume that
  our sextuple is in this subset.} that the six points are in
general position; in particular, this
ensures that these sixty Pascals are pairwise distinct (see~\cite{Pedoe}). 

\subsection{Primal and dual notation}
The notation $L(abcdef)$ for a Pascal line exhibits the
particular hexagon from which it is defined.  In \cite{Ladd}, Christine
Ladd extends this to a {\em primal notation} that gives natural names 
(indexed by
combinations of the points $a, b, c, d, e$ and $f$) for all of the
other nodes and lines in the multimysticum.  The very similar notation
used in \cite{ConwayRyba1} interprets the indices as permutations of the
six points, with the understanding that inverse permutations index the
same object.  For example, the papers~\cite{ConwayRyba1,Ladd} discuss a
\emph{Steiner node} $N(ace)(bdf) = N(eca)(fdb)$ that turns out to
lie on the three Pascal lines $L(abcdef), L(adcfeb)$ and $L(afcbed)$
whose indexing permutations square to $ace.bdf$.  The incidences
between all nodes and lines in the multimysticum are described by 
convenient group theoretic relationships between indexing permutations,
and in particular are invariant under the action of the group
$S_6$ that permutes the six points.

In his construction of the multimysticum Veronese introduces 
\emph{linking lines} and \emph{higher meeting points} to pass between
successive mutations.  In \cite{Ladd}, Ladd observes that there is
no convenient primal notation for these linking objects.  The solution
turns out to be the use of a \emph{dual notation} where we apply
Sylvester's outer automorphism of $S_6$ to the indexing permutations.
The one downside to this notation is that it is now less obvious
which hexagon leads to which Pascal line, but in compensation natural
names become available for Veronese's linking objects.  Of course
these natural names are no longer permutations, since otherwise their
preimages under the outer automorphism would be natural primal names. 

\subsection{The dual notation} 
Define two sets 
\[ \ltr = \{a,b,c,d,e,f\}, \quad \num = \{0,1,2,3,4,5\}. \] 
consisting of six elements each. For any set $X$, let $\SG(X)$ denote the group of bijections
$X \rightarrow X$. The following table\footnote{This table has been reconstructed from the list of six Desargues  
configurations given on~\cite[p.~45]{ConwayRyba2}.}  defines an isomorphism 
$\zeta: \SG(\ltr) \lra \SG(\num)$, 
which is an outer automorphism. For instance, the entry in row $a$ and column
$b$ is to be read as saying that $\zeta$ takes the $2$-cycle
$\transp{a}{b}$ to the element $\cyclett{2}{1}{5}{3}{0}{4}$ of cycle
type $2+2+2$. 

\[ 
\begin{array}{|c|c|c|c|c|c|c|} \hline 
{} & b & c & d & e & f \\ \hline 
a & 21.53.04 & 24.51.03 & 20.54.13 & 25.01.34 & 23.50.14 \\ 
b & {} & 23.54.01 & 25.03.14 & 24.50.13 & 20.51.34 \\ 
c & {} & {} & 21.50.34 & 20.53.14 & 25.04.13 \\ 
d & {} & {} & {} & 23.51.04 & 24.53.01 \\ 
e & {} & {} & {} & {} & 21.54.03 \\ \hline 
\end{array} \label{table.theta} \]  

Now, in order to find the dual notation for an arbitrary Pascal, say $L(acebfd)$, we
find the image of the $6$-cycle $ h = a \, c \, e \, b \, f \, d
\in \SG(\ltr)$. Using the table, it is easy to calculate that 
\[ \zeta(h) = 2\, . \,0 \, 4\, .\, 1 \, 5 \, 3, \quad \text{and} \quad
\zeta(h^{-1}) = 2\, .\, 0 \, 4\, .\, 1 \, 3 \, 5,  \] 
which are elements of cycle type $1+2+3$. In particular, the
$1+2$ part of the image only depends on the
set $\{h,h^{-1}\}$. Hence we relabel $L(acebfd)$ as $\pasc{2}{0}{4}$ or
$\pasc{2}{4}{0}$. In this way, the sixty Pascals are labelled as $\pasc{x}{y}{z}$,
where $x,y,z \in \num$ are distinct elements and the
order of $y,z$ is immaterial. 

\subsection{} 
The action of $\zeta^{-1}$ can be easily worked out from
the same table. For instance, the pair $25$ appears in places $(a,e),
(b,d)$ and $(c,f)$, from which we deduce that $\transp 25 \stackrel{\zeta^{-1}}{\lra} 
\cyclett aebdcf$. In order to find the hexagon
corresponding to an arbitrary Pascal, say $\pasc{3}{1}{5}$, one calculates that 
\[ 3\, .\, 1 \, 5\, . \, 0 \, 2 \, 4 \stackrel{\zeta^{-1}}{\lra} 
a \, d \, f \, c \, e \, b, \] 
and hence $\pasc{3}{1}{5} = L(adfceb)$. Of course, the hexagon is
determined only up to a cyclic shift and reversal. 

\subsection{} \label{section.ordinary.meeting.point} 
There is an alternate way to see the labelling
convention; which will be `generalised' when we define the higher
mutations in Section~\ref{section.multimysticum}. 

The Pascal $\pasc{2}{0}{4} = L(acebfd)$ is the line passing through the 
intersection points $ac \cap bf, ce \cap df, ad \cap be$. We can
introduce `dual notations' for these points as follows: interpret the
first as the element $\transp{a}{c}\, . \, \transp{b}{f} \in \SG(\ltr)$. When
we apply $\zeta$, it gets sent to $\transp{2}{3}\, . \,\transp{0}{4} \in \SG(\num)$.
The same calculation done on the other two
points shows that $\pasc{2}{0}{4}$ is the line passing through $21.04, 23.04$
and $25.04$. The pattern for a general $\pasc{x}{y}{z}$ is now clear;
it passes through $xu.yz, xv.yz, xw.yz$ where $\{u,v,w,x,y,z\} =
\num$. 

Altogether there are $45$ such intersections points $xy.zw$, which
will be called \emph{ordinary} meeting points. Each of them lies on four
Pascals; for instance, $12.35$ lies on $\pasc{1}{3}{5}, \pasc{2}{3}{5}, \pasc{3}{1}{2}$ and
$\pasc{5}{1}{2}$. We will come across \emph{higher} meeting points
later in Section~\ref{section.mutation.eveni}. 

\subsection{} 
We proceed to describe the rest of the elements in $\HM$. The proofs of all the
incidence theorems stated below may be found in
\cite{ConwayRyba1, ConwayRyba2}. In keeping with the conventions of
those papers, the points which participate in $\HM$ will be called `nodes'. The
notations for all the lines begin with $L$, and those for nodes begin
with $N$.  The construction sequence is as follows: 

\begin{equation} \xymatrix{
\text{Pascal lines} \ar[r] \ar[d] & \text{Kirkman nodes} \ar[d] \\ 
\text{Steiner nodes} \ar[d]  & \text{Cayley lines} \ar[d] \\ 
\text{Pl{\"u}cker lines}  & \text{Salmon nodes} }
\label{six.types.HM} \end{equation} 

That is to say, the Kirkman nodes are constructed from Pascal lines using
a concurrency theorem, the Pl{\"u}cker lines are constructed from
Steiner nodes using a collinearity theorem, and so on. 

\subsection{The Kirkman nodes} 
\label{section.kirkman.nodes}

A theorem of Kirkman shows that the 
Pascals $\pasc{0}{1}{2}, \pasc{0}{1}{3}, \pasc{0}{2}{3}$ are concurrent. Their common
point (called a Kirkman node or just a Kirkman) will be denoted by
$\kirk{0}{4}{5}$, where the pair $45$ is obtained by removing $0,1,2,3$
from $\num$. It is understood that a parallel statement is true
for any such indicial pattern; for example, the Pascals 
$\pasc{1}{0}{2}, \pasc{1}{0}{5}, \pasc{1}{2}{5}$ 
are concurrent in the node $\kirk{1}{3}{4}$ etc. Altogether, we have $60$ Kirkman
nodes $\kirk{x}{y}{z}$. 

Each Pascal contains three Kirkmans, and each Kirkman lies on
three Pascals. The incidence pattern follows naturally from the
notation. For instance, $\pasc{1}{2}{3}$ contains the three
Kirkmans $\kirk{1}{0}{4}, \kirk{1}{0}{5}, \kirk{1}{4}{5}$. Similarly, $\kirk{1}{0}{5}$ lies on
$\pasc{1}{2}{3}, \pasc{1}{2}{4}, \pasc{1}{3}{4}$. In particular,
although the Pascals are logically prior, they can in fact be
recovered from the Kirkmans. This dependence will look more
\emph{symmetric} when we come to the higher mutants in the next
section. 

The Pascals and the Kirkmans together form the mutable part
$\MM$ of $\HM$. The fixed part $\FF$ is made of two more kinds of lines
and nodes which are described below. 

\subsection{Steiner nodes and Cayley lines} 
\label{section.SteinerCayley} 
Choose three elements in $\num$, say $0,2,5$. There are three Pascal
lines, namely $\pasc{1}{3}{4}, \pasc{3}{1}{4}, \pasc{4}{1}{3}$, whose labels avoid these
three elements. A theorem of Steiner shows that these three
Pascals are concurrent. Their common point is the Steiner node
$N025$. 

An analogous statement is true of the Kirkman nodes; that is
to say, $\kirk{1}{3}{4}, \kirk{3}{1}{4}, \kirk{4}{1}{3}$ are collinear; they lie
on the Cayley line $L025$.  Thus, there are 
$\binom{6}{3} = 20$ lines and points respectively labelled $Lxyz$ and $Nxyz$, where
$x,y,z \in \num$ are distinct elements whose order is immaterial. 

The Steiner node $Nxyz$ lies on the Cayley line
$Luvw$, whenever $\{x,y,z,u,v,w\} = \num$. For instance, $N024$ lies on $L135$.

\subsection{Pl{\"u}cker lines and Salmon nodes} 
\label{section.PluckerSalmon} 
Choose two elements in $\num$, say $2$ and $4$. This leaves out
$0,1,3,5$. Now consider all Steiner nodes formed using these four
leftover numbers, namely $N013,N015,N035$ and $N135$. A theorem of
Pl{\"u}cker shows these four lie on a line, which we call the
Pl{\"u}cker line $L24$. 

An analogous statement is true of the Cayley lines; that is to
say, $L013,L015,L035$ and $L135$ are concurrent in a point which will be
called the Salmon node $N24$.  There are
$\binom{6}{2} = 15$ such lines and points respectively labelled $Lxy$
and $Nxy$, where $x,y \in \num$ are distinct elements whose order is immaterial. 

This completes the description of $\HM$, with its $60+20+15 = 95$
lines and 95 points. To recapitulate, all the Pascal lines and Kirkman
nodes $\pasc{x}{y}{z}, \kirk{x}{y}{z}$ are in the mutable part $\MM$, whereas all the lines and nodes
defined in Sections~\ref{section.SteinerCayley}--\ref{section.PluckerSalmon} 
are in the fixed part $\FF$. For a \emph{general} choice of six points on the conic,
there are no incidences apart from those already mentioned.

\subsection{} \label{section.classical.references}
The articles by Conway and Ryba~\cite{ConwayRyba1, ConwayRyba2} contain
all of the preceding material in a rather concentrated form. One of the 
best early surveys of the field is due to George Salmon
(see~\cite[Notes]{SalmonConics}). Much of this material is also 
explained by H.~F.~Baker in his note `On the 
\emph{Hexagrammum Mysticum} of Pascal' in~\cite[Note
II]{Baker}. Although his labelling scheme is ostensibly different from 
the one used here; it is based upon the same
foundational notion, namely the outer automorphism of $\SG_6$. Several
more older sources are listed in the bibliographies
of~\cite{ConwayRyba1, ConwayRyba2}. Many of the combinatorial and geometric aspects of the outer
automorphism are discussed\footnote{The outer automorphism of
$\SG_6$ is unique up to inner automorphisms. Hence, all of its
descriptions, wherever they may appear in the literature, are
equivalent to each other up to relabelling.} in~\cite{HMSV}. 

Broadly speaking, there are (at least) two entirely different methods
available for proving all of these incidences. One is to use
Desargues' theorem repeatedly, as in~\cite{ConwayRyba1}
or in the `Notes' by Salmon~\cite{SalmonConics}. 
Another approach, which can be attributed to Cremona
and Richmond, is to use the geometry of lines on a nodal cubic surface 
(see~\cite{Richmond1, Richmond2}). 

The  Steiner-Cayley elements have incidence relations with the
Pascal-Kirkman as well as Pl{\"u}cker-Salmon elements. However, there
are no direct incidence relations going across the latter two; that is to say, 
none of the Pl{\"u}cker lines contains any of the Kirkman nodes and
none of the Salmon nodes lies on any of the Pascal lines. 

Therefore, the Steiner-Cayley elements are in some ways central to the
structure of $\HM$. They will naturally serve as bases for the two ranges
described below. 
\subsection{The Kirkman range} 
\label{section.kirkman.range} 
Consider an arbitrary Kirkman node, say $\kirk{3}{0}{5}$. It lies on the
Cayley line $L124$, which also contains the nodes $N05$ and $N305$. 
Thus we have three points 
\[ N05, \quad N305, \quad \kirk{3}{0}{5}, \] 
on $L124$ which, in effect, create a coordinate system on this line. In
other words, there is a unique isomorphism $f: L124 \lra \PL$ 
such that these points are taken to $\infty, 0, 1$ respectively. Then 
each point $Q \in L124$ is uniquely identified by its image $f(Q) \in
\PL$. 

\begin{figure}
\includegraphics[width=16cm]{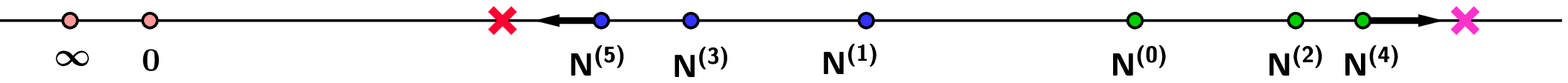}
\caption{The Kirkman range} 
\label{kirkrange.diagram} 
\end{figure}

In the next section we will see that the Kirkman point $\kirk{3}{0}{5}= \hkirk{3}{0}{5}{0}$ has
higher \emph{mutants} $\hkirk{3}{0}{5}{i}$ for $i \geqslant 1$, all of which lie on $L124$. Hence we
have a range 
\[ N05, \quad N305, \quad \kirk{3}{0}{5}, \quad \hkirk{3}{0}{5}{1}, \quad 
\hkirk{3}{0}{5}{2}, \dots \] 
on $L124$. In general, given a Kirkman point $\kirk{x}{y}{z}$,  write 
$\{x,y,z,u,v,w\} = \num$, and consider the range
\begin{equation} 
Nyz, \quad Nxyz, \quad \kirk{x}{y}{z}, \quad \hkirk{x}{y}{z}{1}, \quad
\hkirk{x}{y}{z}{2}, \dots 
\label{kirkman.range} \end{equation}
on the Cayley line $Luvw$. Let $f: Luvw \lra \PL$ denote the isomorphism
which takes the first three points to $\infty, 0, 1$ respectively. Now
the first part of our main theorem is as follows: 

\begin{Theorem}[First Part] \rm 
With notation as above, the range~(\ref{kirkman.range}) is isomorphic
to the Veronese sequence given in~(\ref{Veronese.sequence}). In other words,
we have $f(\hkirk{x}{y}{z}{i}) = \alpha_i$ for $i \geqslant 1$. 
\end{Theorem}  
In particular, as an abstract projective range, the sequence of nodes 
in~(\ref{kirkman.range}) is independent of the indices $x,y,z$, and also
of the initial choice of the sextuple on the conic. 
One can see this schematically\footnote{The diagram is not to scale. Since the point $\infty$ is
  shown, it cannot be.} in Diagram~\ref{kirkrange.diagram}. The
blue nodes $N^{(2r+1)}$ at odd heights have a limiting position (shown 
as a red cross), while the green nodes $N^{(2r)}$ at even heights have
a different limiting position (shown as a pink cross). 

There are a few
variants of this theorem to be described below; a combined proof for all of them will be
given in Section~\ref{section.maintheorem.proof}. 
\subsection{The Pascal range} 
\label{section.pascal.range} 
Instead of a Kirkman node, one can begin with a Pascal
line and carry out exactly the same construction with the roles of
$N$ and $L$ interchanged. Thus, just as in~(\ref{kirkman.range}), we
have a sequence of lines 
\begin{equation} 
Lxyz, \quad Lyz, \quad Lx_{yz}, \quad \hpasc{x}{y}{z}{1}, \quad 
\hpasc{x}{y}{z}{2}, \dots 
\label{pascal.range} \end{equation}
all passing through the Steiner node $Nuvw$. 

Notice that the ranges in~(\ref{kirkman.range})
and~(\ref{pascal.range}) formally look the same from the third element
onwards (with $N$ and $L$ interchanged), but the 
first two elements have been transposed. The Kirkman range starts as $Nyz, Nxyz$, whereas
the Pascal range starts as $Lxyz, Lyz$. This is, of course,
deliberate. 

\begin{Theorem}[Second Part] \rm 
With notation as above, the range~(\ref{pascal.range}) is also 
isomorphic to the Veronese sequence. 
\end{Theorem} 

\begin{figure}
\includegraphics[width=14cm]{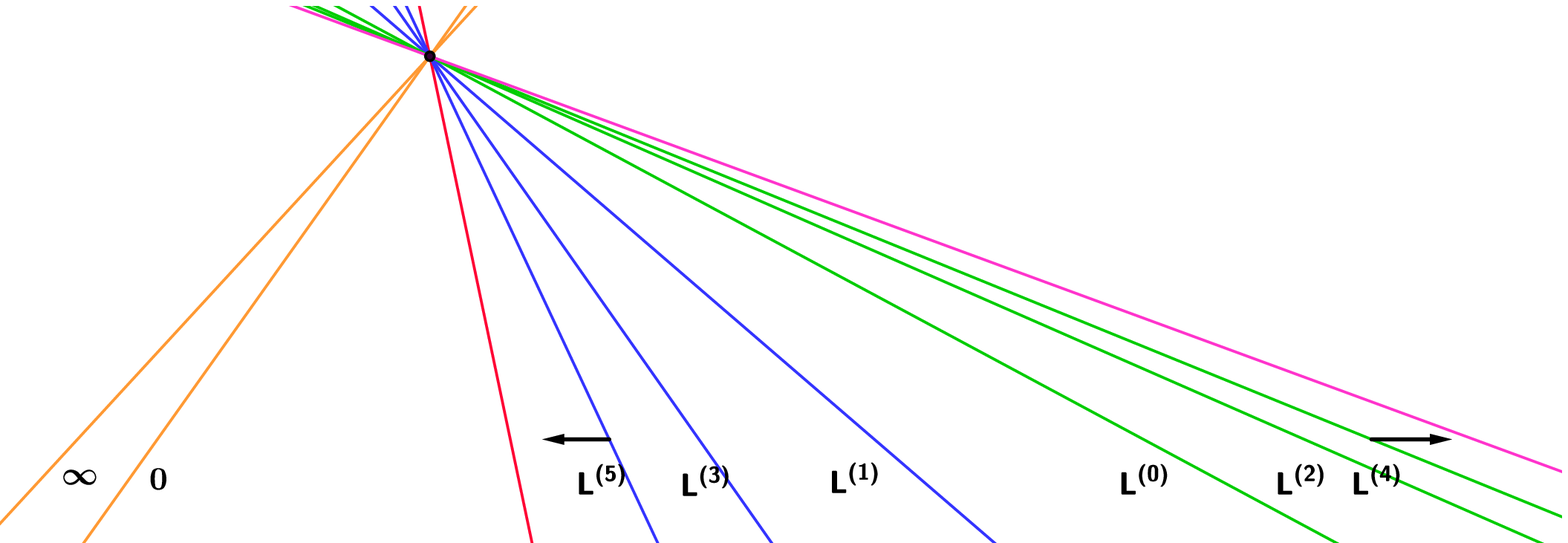}
\caption{The Pascal range} 
\label{pascrange.diagram} 
\end{figure}

The corresponding schema is shown in Diagram~\ref{pascrange.diagram}. The
Pascals at the odd heights (shown in blue) converge to the red line,
whereas those at the even heights (shown in green) converge to the pink
line. 

Either of the ranges above picks up one geometric element from each of the $\MM^{(i)}$. By
contrast, there are two more infinite sequences, namely the `meeting
range' and the `linking range', which are associated to the mutations 
$\MM^{(i)} \lra \MM^{(i+1)}$. As such, their elements could be said
to lie `midway' between two adjacent $\MM^{(i)}$. 
They will also turn out to be isomorphic to the Veronese sequence. 

\section{Veronese mutations and the Multimysticum} 
\label{section.multimysticum} 

In this section, we will describe the process of mutation and in
particular the two new ranges. Once this set-up is in place, the proof of the main theorem
itself will be rather concise. Nearly all the material in this section
may be found in Veronese's original memoir~\cite{Vero}, and also
in the subsequent paper by Ladd~\cite{Ladd}. 

Let $i$ be an index $\geqslant 0$. If $i$ is even, then the
passage from $\MM^{(i)}$ to $\MM^{(i+1)}$ uses an auxiliary system of
`linking lines'; whereas for $i$ odd, it uses a system of `meeting points'. 
Recall that $\MM^{(0)}$ represent the usual Pascal
lines and Kirkman nodes. The higher `Pascals' and `Kirkmans' lie in
$\MM^{(i)}$ for $i \geqslant 1$, but it should be kept in mind that
they do not come from an actual sextuple of points on a conic. 

Along the way, we will refer to several incidence properties of the
new points and lines as they arise. The proofs of all of these may be 
found in~\cite{ConwayRyba2}.

\subsection{} \label{section.mutation.eveni} 
Assume that $i$ is even. The mutation $\MM^{(i)} \lra
\MM^{(i+1)}$ can be schematically represented as follows: 

\begin{equation} \xymatrix{
\text{Kirkman nodes at height $i$} \ar[r] & 
\text{linking lines} \ar[r] & \text{Kirkman nodes at 
  height $i+1$} \ar[d] \\ 
\text{Pascal lines at height $i$} & & 
\quad \text{Pascal lines at height $i+1$}} 
\end{equation} 

To wit, the Kirkmans at height $i$ are used to create a system of linking
lines, which are in turn used to construct Kirkmans at height
$i+1$. These new nodes are then used to construct new Pascals at
height $i+1$. The Pascals at height $i$ play no direct role in the
construction. The recipe is as follows: 
\begin{enumerate} 
\item 
For any four elements in $\num$, say $1,2,3,4$, define the linking
line 
\begin{equation} 
1^{(i+1)} \, 2^{(i+1)} . \, 3^{(i)} \, 4^{(i)} = 
\text{line joining the nodes $\hkirk{3}{1}{2}{i}$ and
  $\hkirk{4}{1}{2}{i}$.} 
\label{linking.line} \end{equation} 
There are $\binom{6}{4} \binom{4}{2} = 90$ such lines. 
\item 
Now define a higher Kirkman node, say $\hkirk{1}{3}{4}{i+1}$, as the common
point of the lines 
\[ 1^{(i+1)} \, x^{(i+1)} . \, 3^{(i)} \, 4^{(i)}, \quad x = 0,2,5. \] 
It is a fact that these three lines are concurrent, and hence this
node is well-defined (see~\cite{ConwayRyba2}). 
\item 
Once all the new Kirkmans are in place, the new Pascals are
constructed from them as in Section~ \ref{section.kirkman.nodes}. That is to
say, define $\hpasc{1}{2}{3}{i+1}$ as the line containing 
\[ \hkirk{1}{0}{4}{i+1}, \quad \hkirk{1}{0}{5}{i+1}, \quad
\hkirk{1}{4}{5}{i+1}. \] 
Again, it is a fact that these three new nodes are collinear, and
hence the new Pascal is well-defined.  
\end{enumerate} 

\begin{figure}
\includegraphics[width=6cm]{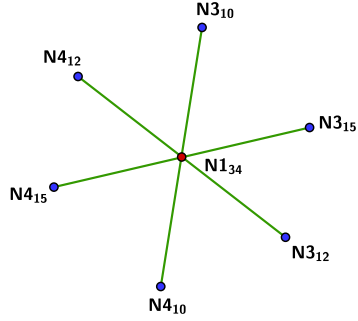}
\caption{For $i$ even, the blue Kirkmans at height $i$ lead to the red 
  Kirkman at height $i+1$.} 
\label{even2oddK.diagram} 
\end{figure}

The first two steps are illustrated in Diagram~\ref{even2oddK.diagram}. The
Kirkmans at height $i$ (shown in blue) are connected by linking lines
(shown in green). Three of these lines are concurrent in the new
Kirkman (shown in red) at height $i+1$. 

\subsection{} \label{section.mutation.oddi} 
Now assume that $i$ is odd. 
The mutation $\MM^{(i)} \lra \MM^{(i+1)}$ is given by a `reciprocal'
procedure. Schematically, we have: 

\begin{equation}  \xymatrix{
\text{Kirkmans at height $i$}  & 
& \text{Kirkmans at height $i+1$} \\ 
\text{Pascals at height $i$} \ar[r] & 
\text{meeting points} \ar[r] & \text{Pascals at height $i+1$}
\ar[u] } \end{equation} 

The Pascals at height $i$ are used to create a system of
meeting points, which are in turn used to construct Pascals at height
$i+1$. These new lines are then used to construct new Kirkmans at
height $i+1$. The Kirkmans at height $i$ play no direct role in the
construction. The recipe, which follows by making cosmetic changes to
the earlier recipe, is as follows: 

\begin{enumerate} 
\item 
For any four elements in $\num$, say $1,2,3,4$, construct the higher meeting
point 
\begin{equation} 
1^{(i+1)} \, 2^{(i+1)} . \, 3^{(i)} \, 4^{(i)} = 
\text{point of intersection of the lines $\hpasc{3}{1}{2}{i}$ and
  $\hpasc{4}{1}{2}{i}$.} 
\label{meeting.point} \end{equation} 
There are $90$ such points. They are notationally distinguished from the
ordinary meeting points in Section~\ref{section.ordinary.meeting.point}, and hence we will omit
the adjective `higher' if no confusion is likely. 
\item 
Now define a higher Pascal, say $\hpasc{1}{3}{4}{i+1}$, as the line passing through
the meeting points 
\[ 1^{(i+1)} \, x^{(i+1)} . \, 3^{(i)} \, 4^{(i)}, \quad x = 0,2,5. \] 
These three points are in fact collinear, hence the Pascal is
well-defined. 
\item 
Once the new Pascals are in place, the new Kirkmans are
constructed from them as in Section~ \ref{section.kirkman.nodes}. That is to
say, define $\hkirk 123{i+1}$ as the common point of the concurrent lines 
\[ \hpasc{1}{0}{4}{i+1}, \quad \hpasc{1}{0}{5}{i+1}, \quad
\hpasc{1}{4}{5}{i+1}. \] 
\end{enumerate} 
\begin{figure}
\includegraphics[width=8cm]{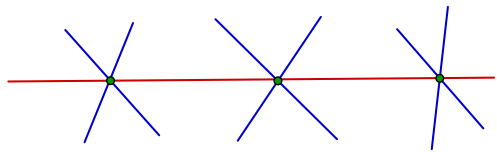}
\caption{For $i$ odd, the blue Pascals at height $i$ lead to the red 
 Pascal at height $i+1$.} 
\label{odd2evenP.diagram} 
\end{figure} 

The first two steps are shown in Diagram~\ref{odd2evenP.diagram}. The
Pascals at height $i$ (shown in blue) intersect in meeting points 
(shown in green). Three of the meeting points are collinear in the new
Pascal (shown in red) at height $i+1$. We have omitted all the labels, 
since they can be obtained from Diagram~\ref{even2oddK.diagram} by replacing each $N$ 
with an $L$. 

This completes the process of defining the mutations. It is a fact
that all the new Kirkmans and Pascals have exactly the same incidence relations
with the Cayley lines and Steiner nodes respectively. For
instance, the Kirkman $\hkirk{1}{3}{4}{i}$ lies on the Cayley line $L025$
for all values of $i \geqslant 0$. Similarly, the Pascal
$\hpasc{1}{3}{4}{i}$ passes
through the Steiner node $N025$ for all values of $i \geqslant 0$. This ensures
that the Kirkman range in Section~\ref{section.kirkman.range}, and the Pascal
range in Section~\ref{section.pascal.range} are well-defined. 

\subsection{} Notice that the linking lines in~(\ref{linking.line}) and 
meeting points in~(\ref{meeting.point}) have formally the same
notation. It is the parity of $i$ which tells them apart. Moreover,
the symmetry $xy.zw = zw.xy$ which holds for ordinary meeting points in
Section~\ref{section.ordinary.meeting.point} is no longer valid for
higher meeting points and linking lines; that is to say, 
\[ 1^{(i+1)} \, 2^{(i+1)} . \, 3^{(i)} \, 4^{(i)} \neq 3^{(i+1)} \,
4^{(i+1)} . \, 1^{(i)} \, 2^{(i)}. \] 
Thus there are $90$ higher meeting points for each odd $i$ and $90$ linking lines
for each even $i$, but
only $45$ ordinary meeting points.   

The larger number of higher meeting 
points proves that no set of $60$ higher Pascal lines can arise as the
Pascal lines corresponding to six points on any conic.  (Otherwise, the
$45$ meeting points corresponding to these Pascal lines would have to match
the $90$ higher meeting points obtained as intersections of higher Pascal lines.
This is clearly impossible since $45 \neq 90$.)  A similar argument
shows that no set of Kirkman nodes (or higher Kirkman nodes) can arise
as the $60$ points obtained by applying Brianchon's theorem 
to a set of six tangents to a conic.

We now define the Veronese nodes and Ladd lines. All the incidences referred
to are proved in~\cite{ConwayRyba2}. 

\subsection{The Ladd lines} 
Assume $i$ to be odd, and consider the two meeting points 
\[ 1^{(i+1)} \, 2^{(i+1)} . \, 3^{(i)} \, 4^{(i)} \quad \text{and}
\quad 3^{(i+1)} \, 4^{(i+1)} . \, 1^{(i)} \, 2^{(i)}. \] 
It is a remarkable fact that the line joining them is \emph{independent of} $i$. It is called the
Ladd line $L12.34$. There are $45$ such lines $Lxy.zw$. 

The Ladd line $L12.34$ contains the Salmon node $N05$ and the ordinary meeting point
$12.34$. Define the Pl{\"u}cker-Ladd node $\PLN{1}{2}{3}{4}$ to be the
intersection of the Pl{\"u}cker line $L34$ and the Ladd line
$L12.34$. At this point, we have all the ingredients necessary to create a range
on $L12.34$. Write 
\[ p_i = 3^{(i+1)} \, 4^{(i+1)} . \, 1^{(i)} \, 2^{(i)} \quad
\text{and} \quad 
q_i = 1^{(i+1)} \, 2^{(i+1)} . \, 3^{(i)} \, 4^{(i)}, \] 
and consider the range 
\begin{equation} 
N05, \quad \PLN{1}{2}{3}{4}, \quad 12.34, \quad p_1, \quad q_1,
\quad p_3, \quad q_3, \quad p_5, \quad q_5, \dots 
\label{meeting.range} \end{equation} 
which will be called the `meeting range', since after the initial
stretch of two points it is entirely made of meeting points. Notice
that there are $90$ such ranges, although there are only $45$ Ladd
lines. This is so because $\PLN{1}{2}{3}{4}$ and $\PLN{3}{4}{1}{2}$
are distinct points belonging to distinct ranges on the same
line. This circumstance will prove useful in the proof of the main theorem
(see Section~\ref{section.proof.step2} below). 

\subsection{The Veronese nodes} 
Assume $i$ to be even, and consider the two linking lines 
\[ 1^{(i+1)} \, 2^{(i+1)} . \, 3^{(i)} \, 4^{(i)} \quad \text{and}
\quad 3^{(i+1)} \, 4^{(i+1)} . \, 1^{(i)} \, 2^{(i)}. \] 
Their point of intersection is independent of $i$. It is called the
Veronese node $N12.34$. There are $45$ such nodes $Nxy.zw$. 

The Pl{\"u}cker line $L05$ passes through the Veronese node
$N12.34$. Define the Salmon-Veronese line $\SVL{1}{2}{3}{4}$ to be the
line joining the Salmon node $N34$ and the Veronese node $N12.34$. 
Now write 
\[ r_i = 3^{(i+1)} \, 4^{(i+1)} . \, 1^{(i)} \, 2^{(i)} \quad
\text{and} \quad s_i = 1^{(i+1)} \, 2^{(i+1)} . \, 3^{(i)} \, 4^{(i)}, \] 
and consider the range of lines 
\begin{equation} 
\SVL{1}{2}{3}{4}, \quad L05, \quad r_0, \quad s_0, \quad r_2, \quad
s_2, \quad r_4, \quad s_4, \dots 
\label{linking.range} \end{equation} 
all passing through $N12.34$. In analogy with the above, this will be
called a `linking range'. There are $90$ such ranges with two of them
based upon each Veronese node. 

The constructions in~(\ref{meeting.range}) and~(\ref{linking.range})
are not perfectly parallel. Apart from the transposition in the first
two elements, there is an additional asymmetry which comes from the fact
that there is no linking line on an equal footing with the meeting point
$12.34$. 

\begin{Theorem}[Third Part] \rm 
With notation as above, the ranges in~(\ref{meeting.range})
and~(\ref{linking.range}) are also isomorphic to the Veronese
sequence. 
\end{Theorem} 

Thus, we have altogether $60+90=150$ point ranges, and the same number
of line ranges, all isomorphic to the Veronese sequence. Although we
did not initially say so, the Ladd lines and Veronese nodes can also
been seen as belonging to the fixed part $\FF$. 

\section{Proof of the Main Theorem} 
\label{section.maintheorem.proof} 

\subsection{}  In this section, we will finally prove the main theorem. 
It will be convenient to label the four ranges. Write $\{u,v,w,x,y,z\}
= \num$, and let 
\begin{equation} 
\begin{aligned} 
\KKK{x}{y}{z} & = \text{Kirkman range in~(\ref{kirkman.range}) based upon the 
  Cayley line $Luvw$,} \\ 
\PPP{x}{y}{z} & = \text{Pascal range in~(\ref{pascal.range}) based upon the 
  Steiner node $Nuvw$,} \\ 
\MMM{x}{y}{z}{w} & = \text{Meeting range in~(\ref{meeting.range}) based upon the 
  Ladd line $Lxy.zw$,} \\ 
\LLL{x}{y}{z}{w} & = \text{Linking range in~(\ref{linking.range}) based upon the 
  Veronese node $Nxy.zw$.} 
\end{aligned} \label{four.ranges} \end{equation}  
Recall that we get altogether $150$ point ranges and the 
same number of line ranges by assigning the letters $u, v, \dots, z$ to 
$0,1, \dots, 5$ in all possible ways. 

Now the main theorem will proved in the following steps: 
\begin{enumerate} 
\item First, we show that all the $300$ ranges in~(\ref{four.ranges})
  are isomorphic to each other. This implies that any of them is isomorphic to a range 
\begin{equation} 
\infty, \quad 0, \quad \beta_0=1, \quad \beta_1, \quad
\beta_2, \quad \beta_3, \dots 
\label{beta.range} \end{equation} 
on $\PL$ for some as yet unknown complex numbers $\beta_i$. 
\item
Next, we show that: 
\begin{itemize} 
\item 
For $i$ odd, we have $\beta_i+ \beta_{i+1} = 2$. 
\item 
As $m$ ranges over the nonnegative integers, the sum $1/\beta_{2m} +
1/\beta_{2m+1}$ remains independent of $m$. 
\end{itemize} 
\item Finally, we show that $\beta_1 = 1/2$, and hence this sum is in
 fact $3$. This shows that~(\ref{beta.range}) is identical to
 the Veronese sequence in (\ref{Veronese.sequence}), which completes the argument. 
\end{enumerate} 

\subsection{Aligned ranges} \label{section.aligned.ranges} 
\begin{figure}
\includegraphics[width=8cm]{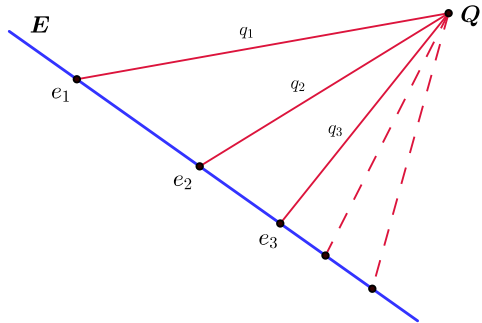}
\caption{Aligned (and hence isomorphic) ranges} 
\label{incident_ranges.diagram} 
\end{figure}

Let $Q$ and $E$ be respectively a point and a line in
the projective plane, such that $Q \notin E$. Assume that $(q_1,q_2, \dots)$
and $(e_1, e_2, \dots)$ are ranges based upon $Q$ and $E$ respectively; that
is to say, $q_i$ are lines through $Q$ and $e_i$ are points on
$E$. If $e_i$ is incident with $q_i$ for all $i$, then the ranges are said to be
aligned (see Diagram~\ref{incident_ranges.diagram}). In particular, they are then isomorphic. We will
complete \textbf{step one} by using this technique repeatedly. 

Many of the proofs below depend upon exploiting a specific 
\emph{indicial pattern}. In such cases, for the sake of vividness, we will use
specific numbers from $0$ through $5$ instead of letters $u,v, \dots,
z$. This should make the proofs easier to follow. 
We will use the incidence
properties of points and lines from
Sections~\ref{section.kirkman.nodes}--\ref{section.PluckerSalmon}, and
also those involving higher mutants from
Sections~\ref{section.mutation.eveni}--\ref{section.mutation.oddi}. 

\subsection{} \label{section.aligned.example} 
Notice that the Kirkman $\kirk{2}{0}{4}$ lies
on the Pascal $\pasc{2}{1}{5}$. This incidence lifts to an alignment
between the corresponding Kirkman and Pascal ranges:
$\KKK{2}{0}{4}$ and
$\PPP{2}{1}{5}$. Indeed, the two ranges respectively correspond to the two
rows: 
\[ \begin{array}{|c|c|c|c|c} \hline 
N04 & N204 & \kirk{2}{0}{4} & \hkirk{2}{0}{4}{1} & \dots \\ \hline 
L215 & L15 & \pasc{2}{1}{5} & \hpasc{2}{1}{5}{2} & \dots \\ \hline 
\end{array} \] 
Now it is merely a matter of checking that the node and the line in
any column are incident. For instance, the 
node $N204$ in the second column is incident with the line $L15$ below
it, because $0,2,4$ are disjoint from $1,5$. In summary, if we take
$Q = N034$ (the base of the Pascal range), $E = L125$ (the base of the
Kirkman range) and think of the rows as $(e_1, e_2, \dots), (q_1, q_2,
\dots)$ respectively, then everything matches
Diagram~\ref{incident_ranges.diagram}. 

This observation, used repeatedly, will allow us to correlate a large
number of ranges. 
\begin{Lemma} \rm 
Fix an element $x \in \num$. Then all Kirkman ranges
$\KKK{x}{\bullet}{\bullet}$ and Pascal ranges $\PPP{x}{\bullet}{\bullet}$ are
isomorphic.
\end{Lemma} 
In either instance, the bullets stand for any two indices different from
$x$. 

\Proof 
Consider the following chain of alternating Pascals and Kirkmans: 
\[ \pasc{0}{1}{2}, \quad \kirk{0}{3}{4}, \quad \pasc{0}{5}{1}, \quad 
\kirk{0}{2}{3}, \quad \pasc{0}{4}{5}. \] 
It has the property that any adjacent node-line pair is
incident; in fact the chain is constructed by keeping $x=0$ fixed, and going
through the sequence $1-2-3-4-5-1$ cyclically in disjoint pairs. Hence, the corresponding ranges
\[ \PPP{0}{1}{2}, \quad \KKK{0}{3}{4}, \quad \PPP{0}{5}{1}, \quad 
\KKK{0}{2}{3}, \quad \PPP{0}{4}{5}, \]  
are isomorphic. Now observe that this captures all possible
indicial patterns. For instance, two Pascal ranges
$\PPP{x}{\bullet}{\bullet}$ must either look like $\PPP{0}{1}{2}$ and $\PPP{0}{5}{1}$ with one overlap
in the remaining indices, or like $\PPP{0}{1}{2}$ and $\PPP{0}{4}{5}$ with
no such overlap. We have shown they are isomorphic in either case. Hence all 
Kirkman ranges $\KKK{0}{\bullet}{\bullet}$ are isomorphic to this common
Pascal range, and also to each other. \qed 

We can picture the Pascal and
Kirkman ranges (altogether $60$ and $60$) as being distributed across six
islands, corresponding to the values $x=0,1,\dots,5$. By the lemma above, we know
that the ranges on each island ($10$ and $10$) are isomorphic amongst themselves. 
The meeting and linking ranges will allow us to pass from one island to another. 

\begin{Lemma} \rm 
Let $x,y,z,w \in \num$. Then 
\begin{enumerate} 
\item the Kirkman range $\KKK{x}{y}{z}$ is aligned with the linking
  range $\LLL{x}{w}{y}{z}$; and 
\item the Pascal range $\PPP{x}{y}{z}$ is aligned with the meeting
  range $\MMM{x}{w}{y}{z}$. 
\end{enumerate} \end{Lemma} 
\Proof The argument is very similar to the one in
Section~\ref{section.aligned.example}, and merely amounts to
checking that the corresponding lines and points are incident. We will
give an illustration for (1). Let $x,y,z,w = 1,3,4,2$ respectively. 
Then the fifth node in $\KKK{1}{3}{4}$ and the
fifth line in $\LLL{1}{2}{3}{4}$ are respectively 
\[ \hkirk{1}{3}{4}{2} \quad \text{and} \quad 3^{(3)} 4^{(3)}. 1^{(2)}
2^{(2)}, \] 
which are incident by the very definition of the linking line. The
remaining verifications are equally routine, and we leave them to the
reader. \qed 

One can now pass from one island to another. 
\begin{Lemma} \rm 
Let $x,y,z,w \in \num$. Then 
\begin{enumerate} 
\item 
the Kirkman ranges $\KKK{x}{y}{z}$ and $\KKK{w}{y}{z}$ are isomorphic,
and 
\item 
the Pascal ranges $\PPP{x}{y}{z}$ and $\PPP{w}{y}{z}$ are isomorphic. 
\end{enumerate} \end{Lemma} 

\Proof The two Kirkman ranges are isomorphic because each is aligned
with $\LLL{x}{w}{y}{z}$. The argument for Pascal ranges is essentially
the same. \qed 

The preceding lemmas culminate in the following proposition,
which completes \textbf{step one} in the proof of the main theorem.   
\begin{Proposition} \rm 
All the $300$ ranges in~(\ref{four.ranges}) are isomorphic to each other. 
\qed \end{Proposition} 

\subsection{} \label{section.proof.step2}
Each of the expressions involved in
\textbf{step two} involves two adjacent values of $\beta_i$. The trick
is to use an involution (i.e., an automorphism of order $2$) which
will interchange $p_i$ with $q_i$ in the meeting range, and similarly
$r_i$ with $s_i$ in the linking range. 

Recall that each automorphism of $\PL$ is given by a fractional linear
transformation 
\[ \PL \stackrel{f}{\lra} \PL, \qquad f(z) = \frac{pz+q}{rz+s}. \] 
The matrix
$\left[ \begin{array}{cc} p & q \\ r & s \end{array} \right]$ is
nonsingular, and determined up to a nonzero scalar. 

\begin{Proposition} \rm 
We have $\beta_i + \beta_{i+1} = 2$ for all
odd values of $i$. 
\end{Proposition} 
\Proof 
Consider the isomorphic ranges $\MMM{1}{2}{3}{4}$ and
$\MMM{3}{4}{1}{2}$, both of which are based upon the Ladd line $L12.34$. They
are respectively shown in the rows below: 
\[ \begin{array}{|c|c|c|c|c|c|c|c|} \hline 
N05 & N12.3'4' & 12.34 & p_1 & q_1 & p_2 & q_2 & \dots \\ \hline 
N05 & N34.1'2' & 34.12 & q_1 & p_1 & q_2 & p_2 & \dots \\ \hline 
\end{array} \] 
Observe that the first and the third nodes coincide, and fourth
onwards they get interchanged in pairs. 
Now fix coordinates on $L12.34$ such that the top row gets identified
with the sequence~(\ref{beta.range}). Then the rows appear as 
\[ \begin{array}{|c|c|c|c|c|c|c|c|} \hline 
\infty & 0 & 1 & \beta_1 & \beta_2 & \beta_3 & \beta_4 & \dots \\
     \hline 
\infty & \mu & 1 & \beta_2 & \beta_1 & \beta_4 & \beta_3 & \dots \\ \hline 
\end{array} \] 
for some constant $\mu$. Hence there exists a fractional linear
transformation $f(z) = \frac{pz+q}{rz+s}$ which takes the first row to
the second. Now $f(\infty) = \infty$ implies that $r=0$, and we
may assume $s=1$. Then $f(1) = 1$ implies that $f(z) = pz + (1-p)$ for
some $p$. Since $f(f(z)) = z$ and $f$ is not the identity, we must have $p
= -1$, i.e., $f(z) = 2-z$. Hence it follows that $\beta_{i+1} =
f(\beta_i) = 2- \beta_i$ for all odd values of $i$. \qed 

Now we use a similar argument on the linking range. 

\begin{Proposition} \rm 
As $m$ ranges over the nonnegative integers, the sum $1/\beta_{2m} +
1/\beta_{2m+1}$ remains independent of $m$. 
\end{Proposition} 

\Proof 
Consider the isomorphic ranges $\LLL{1}{2}{3}{4}$ and
$\LLL{3}{4}{1}{2}$, both of which are based upon the Veronese node $N12.34$. They
are respectively shown in the rows below: 
\[ \begin{array}{|c|c|c|c|c|c|c|} \hline 
L12.3'4' & L05 & r_0 & s_0 & r_2 & s_2 & \dots \\ \hline 
L34.1'2' & L05 & s_0 & r_0 & s_2 & r_2 & \dots \\ \hline 
\end{array} \] 
The second entry is the same in both ranges, and from the third
onwards they get interchanged in pairs. Now fix coordinates on the
planar pencil of lines through $N12.34$ so that the top row is identified
with the sequence~(\ref{beta.range}). Then the rows appear as 
\[ \begin{array}{|c|c|c|c|c|c|c|} \hline 
\infty & 0 & 1 & \beta_1 & \beta_2 & \beta_3 & \dots \\
     \hline 
\mu & 0 & \beta_1 & 1 & \beta_3 & \beta_2 & \dots \\ \hline 
\end{array} \] 
for some constant $\mu$. As above, there must exists a fractional linear
transformation $f(z) = \frac{pz+q}{rz+s}$ which takes the first row to
the second. Now $f(0) =0$ implies that $q=0$ and we may assume
$p=1$. Since $f$ is a non-identity function such that $f(f(z)) = z$,
we must have $s=-1$ and hence $f(z) = z/(rz-1)$ for some $r$. Now 
$f(\beta_2) = \beta_3$ implies that $r = 1/\beta_2 + 1/\beta_3$. The
same argument applies to all such pairs, and thus we get 
\[ r = 1 + \frac{1}{\beta_1} = \frac{1}{\beta_2} + \frac{1}{\beta_3} =
\frac{1}{\beta_4} + \frac{1}{\beta_5} = \dots \] 
\qed 

\subsection{} 
It only remains to find $\beta_1$. Since the anticipated answer is
$1/2$, the reader may have already guessed that 
\emph{harmonicity} is implicated. Consider the first four elements in
each of the Kirkman ranges $\KKK{0}{2}{3}$ and $\KKK{1}{2}{3}$, namely 
\begin{equation} 
N23, \quad N023, \quad \kirk{0}{2}{3}, \quad \hkirk{0}{2}{3}{1} \qquad
\text{and} \qquad N23, \quad N123, \quad \kirk{1}{2}{3}, \quad
\hkirk{1}{2}{3}{1}. 
\label{two.quadruples} \end{equation} 
We may identify both of them with the range $(\infty,0,1,\beta_1)$. We will construct a 
nontrivial automorphism of this quadruple. Join the four nodes of the
first quadruple to the (ordinary) meeting point $01.45$. This creates
the line quadruple: 
\begin{equation} 
L01.45, \quad \pasc{1}{4}{5}, \quad \pasc{0}{4}{5}, \quad
0^{(1)}1^{(1)}.\, 2^{(0)} 3^{(0)} 
\label{harmonic.quadruple} \end{equation} 
(for the last incidence, see the last line on `Elevation' on page 48 of 
\cite{ConwayRyba2}). 

Now observe that if we now join the second quadruple in
(\ref{two.quadruples}) to the same meeting point, then the middle two elements in~(\ref{harmonic.quadruple}) get
interchanged and those at either end remain unchanged. It follows
that $(\infty,0,1,\beta_1)$ is isomorphic to
$(\infty,1,0,\beta_1)$; that is to say, this is a harmonic quadruple. The only possible linear fractional
transformation which takes the first to the second is $f(z) =
1-z$. Hence $f(\beta_1) = 1-\beta_1 = \beta_1$, which forces 
$\beta_1 = 1/2$. This completes the proof of the main theorem. \qed 

\subsection{} Although our theorem shows that many of the infinite
ranges that belong to the multimysticum are absolutely invariant, this
is certainly not the case for all ranges within the system.  
For example, the Kirkman ranges $\KKK 012, \KKK 102, \KKK 201$
all lie on the Cayley line $L345$, but there is no fixed set of 
projective coordinates for the union of these ranges.
Moreover, there are other points on this Cayley line such as the
anti-Kirkman nodes $\pP {\oo{N}}012$, $\pP {\oo{N}}102$ and $\pP {\oo{N}}201$ 
(see \cite[p.~48]{ConwayRyba2}) which cannot be added to any of these ranges
without sacrificing absolute invariance.

\subsection{} 
The theorem opens up several new lines of inquiry. For instance, it is
clear that some of the terms of the Veronese sequence become zero or
undefined if the base field is of positive characteristic. It would be
of great interest to study the structure of the multimysticum in such
cases. We hope to pick up this thread in a possible sequel to this
paper. 

This paper is based upon a rather more compact version written jointly
by the two authors with Professor John Conway. We should like to thank
him for several helpful discussions. 

{}

{
    \bigskip
      \footnotesize
\noindent
Jaydeep Chipalkatti,
\textsc{Department of Mathematics, Machray Hall, University of Manitoba,
Winnipeg, MB R3T 2N2, Canada}\\ 
{E-mail address}: \texttt{jaydeep.chipalkatti@umanitoba.ca}

\bigskip 

\noindent
Alex Ryba,
\textsc{Department of Computer Science, Queens College, CUNY, 65-30
  Kissena Boulevard, Flushing, NY 11367, USA} \\ 
    {E-mail address}: \texttt{ryba@cs.qc.cuny.edu}
 }

\bigskip 

\centerline{--} 


\begin{thebibliography}{99}

\bibitem{Baker} \mbox{H. F. Baker},
\emph{Principles of Geometry}, volume 2,
University Press, Cambridge, 1922.

\bibitem{ConwayRyba1}
\mbox{J. H. Conway and A. J. E. Ryba},
The Pascal Mysticum Demystified. 
{\em The Mathematical Intelligencer,\ }
{\bf 34}, no.\ 3 (2012), 4--8.

\bibitem{ConwayRyba2}
\mbox{J. H. Conway and A. J. E. Ryba},
Extending the Pascal Mysticum. 
{\em The Mathematical Intelligencer,\ }
{\bf 35}, no.\ 2 (2013), 44--51.

\bibitem{Coxeter} \mbox{H. S. M. Coxeter},
\emph{The real projective plane}, McGrawHill, New York, 1949.

\bibitem{HMSV} 
\mbox{B.~Howard, J.~Millson, A.~Snowden, and R.~Vakil}, 
A description of the outer automorphism of $S_6$, and the
 invariants of six points in projective space.
{\em J.~Combin.~Theory Ser.~A}, {\bf 115}, no.~7 (2008),
1296--1303. 

\bibitem{Ladd} \mbox{C. Ladd},
The Pascal Hexagram, {\em Amer. J. Math.,\ } {\bf 2}  (1879), 1--12. 

\bibitem{Pedoe} \mbox{D.~Pedoe}, 
How many Pascal lines has a sixpoint?
{\em The Mathematical Gazette,\ }, {\bf 25} (1941), 110--111. 

\bibitem{Richmond1} 
\mbox{H.~W.~Richmond}, 
\newblock A symmetrical system of equations 
of the lines on a cubic surface which has a conical point. 
\newblock {\em The Quarterly Journal of Pure and Applied 
Mathematics}, {\bf 23} (1889), 170--179. 

\bibitem{Richmond2} 
\mbox{H.~W.~Richmond},  
\newblock The figure formed from six points in space of four dimensions. 
\newblock {\em Math.~Annalen}, {\bf LIII} (1869), 161--176. 

\bibitem{SalmonConics} \mbox{G.~Salmon}, 
\emph{A Treatise on Conic Sections}, 
Reprint of the 6th ed.~by Chelsea Publishing Co., New York, 2005. 

\bibitem{Vero} \mbox{G. Veronese},
Nuovi teoremi sull'Hexagrammum mysticum, 
{\em Memorie della Reale Accademia dei Lincei,\ } 
{\bf 3}, no.~1 (1877), 649--703. 

\vspace{1cm} 

\end{thebibliography}
\end{document}